\documentclass[11pt]{article}
\usepackage{amsmath}
\usepackage{amssymb}
\usepackage{epsfig}

\newcommand{\Ch}{\mathit{Ch}}                   
\newcommand{\CS}{\mathit{Ch_{QS^0}}}            

\newcommand{\A}{\mathcal{A}}                    
\newcommand{\An}[1][n]{\A^{(#1)}}               

\newcommand{\var}{\displaystyle\mathop{{\rm var}}}   
\newcommand{\colim}{\displaystyle\mathop{\rm colim}} 
\newcommand{\dlim}{\displaystyle\mathop{\rm lim}}    
\newcommand{\Vark}[1]{\var(#1)}                      
\newcommand{\tensor}{\displaystyle\mathop{\otimes}}  

\newcommand{\lra}{\longrightarrow}
\newcommand{\bra}[1]{\buildrel{#1}\over\longrightarrow}
\newcommand{\Ehat}{\widehat{E(n)}}
\newcommand{\E}{\underline{E\!}\,}
\newcommand{\FF}{\underline{F\!}\,}
\newcommand{\EM}{\underline{H\!}\,}
\newcommand{\qed}{\hfill$\square$}
\newcommand{\Image}{\mbox{Im}}
\newcommand{\F}{\mathbb{F}_p}
\newcommand{\Z}{\mathbb{Z}}

\newcommand{\onto}{\twoheadrightarrow}

\newtheorem{thm}{Theorem}[section]
\newtheorem{cor}[thm]{Corollary}
\newtheorem{lem}[thm]{Lemma}
\newtheorem{prop}[thm]{Proposition}
\newtheorem{example}[thm]{Example}
\newtheorem{defn}[thm]{Definition}

\newtheorem{rem}[thm]{Remark}
\newtheorem{warning}[thm]{Warning}

\title{\Huge Subalgebras of group cohomology\\
defined by infinite loop spaces}

\author{John Robert Hunton \& Bj\"orn Schuster}

\date{$30^{\rm th}$ November 2001}
\begin{document}

\maketitle

\begin{abstract}
We study natural subalgebras $Ch_E(G)$ of group cohomology defined in terms
of infinite loop spaces $E$ and give representation theoretic descriptions of
those based on $QS^0$ and the Johnson-Wilson theories $E(n)$. We describe the
subalgebras arising from the Brown-Peterson spectra $BP$ and as a result give
a simple reproof of Yagita's theorem that the image of $BP^*(BG)$ in
$H^*(BG;\F)$ is $F$-isomorphic to the whole cohomology ring; the same result
is shown to hold with $BP$ replaced by any complex oriented theory $E$ with a
map of ring spectra $E\rightarrow H\F$ which is non-trivial in homotopy. We
also extend the constructions to define subalgebras of $H^*(X;\F)$ for any
space $X$; when $X$ is finite we show that the subalgebras $Ch_{E(n)}(X)$
give a natural unstable chromatic filtration of $H^*(X;\F)$.
\end{abstract}

\bibliographystyle{unsrt}
\section{Introduction}

The main aim of this paper is to study certain natural subalgebras $Ch_E(G)$
of $H^*(G)$, the cohomology of a finite group $G$ with coefficients in $\F$,
which are defined using infinite loop spaces associated to a spectrum $E$.
Though not originally defined in this manner, perhaps the first example of
such a subalgebra is that of the {\em Chern subring}, $Ch(G)$, defined and
studied in detail by Thomas \cite{Thomas} and others.
The Chern subring is the subring of $H^*(G)$ generated by the Chern classes
of its irreducible complex representations; it is in some sense that part of
group cohomology most accessible using the tools of representation theory.

In \cite{CCC} the construction $Ch_E(G)$ was introduced for a group $G$
and a spectrum $E$ and was shown to give a subalgebra of $H^*(G)$, natural
in $G$ and closed under the action of the Steenrod algebra. We review
this construction and develop further some of its basic properties in
Section~\ref{basic}. Taking $E$ to be the spectrum  representing complex
$K$-theory,  $Ch_K(G)$ was shown in \cite{CCC} to coincide with Thomas'
Chern  subring  $Ch(G)$. The main goal of \cite{CCC} however was to study the
subalgebras  $Ch_{\Ehat}(G)$ based on Baker and W\"urgler's spectra $\Ehat$,
$n=1,2,\ldots\,$, see \cite{bwa,bwb}, and to give these subalgebras a
representation theoretic description as well. The main tools used were
a \lq Hopf ring\rq\ analysis of the spaces in the $\Omega$-spectra for
$\Ehat$ \cite{hun,HT2,wkn}, the generalised character theory describing
$(\Ehat)^*(BG)$ of \cite{HKR} and the work of Green and Leary \cite{VCh}
on the varieties associated to subrings of $H^*(G)$. The reinterpretation
of $Ch(G)$ as $Ch_K(G)$ was similar, though easier, using the well known
descriptions of the cohomology of $BU$ and Atiyah's theorem on $K^*(BG)$
\cite{Atiyah}.

Despite the work of \cite{CCC}, a number of basic aspects of the
subalgebras $Ch_E(G)$ remained mysterious. In particular, it was not
clear exactly how the subalgebras were related when the spectrum
$E$ was changed, even to a closely related spectrum, nor was it clear what
subalgebras were produced when some very basic examples of spectra were
taken, such as the Brown-Peterson spectra $BP$, the (uncomplete)
Johnson-Wilson theories $E(n)$ \cite{JW} or the universal example $QS^0$.

When a spectrum $E$ has an associated theorem describing $E^*(BG)$ as a
functor from groups to algebras, just as Atiyah's theorem \cite{Atiyah}
identifies $K^*(BG)$ with a certain completion $\widehat{R(G)}$ of the
complex representation ring, it is reasonable to suppose there is a
related theorem describing $Ch_E(G)$, provided a good enough understanding
of the spaces in the $E$-theory $\Omega$-spectrum can be established.
This was the case for the $\Ehat$ spectra in \cite{CCC} using \cite{HKR}
and \cite{hun}.
In a similar way, we show below that a description of $Ch_{QS^0}(G)$ can
be made using the Segal conjecture \cite{carlsson} and classical work on
the homology of $QS^0$ \cite{Priddy}.
Perhaps of more interest though are the cases of subalgebras based on
$BP$ or $E(n)$ as these spectra have at present no associated character
theory or analogue of Atiyah's, Carlsson's or Hopkins-Kuhn-Ravenel's
theorems \cite{Atiyah,carlsson,HKR}.

The case of $BP$-theory turns out to yield results through one further
property of the $Ch_E(G)$ construction, not considered in \cite{CCC}.
By its nature, $Ch_E(G)$ is essentially a construction in {\em unstable\/}
homotopy, and as such there is an associated notion of {\em stabilisation}.
Using this idea we show that $Ch_{BP}(G)$ is in fact $F$-isomorphic to
$H^*(BG)$. Moreover, the same is true for $Ch_E(G)$ whenever $E$ is a
complex oriented spectrum with a map of ring spectra
$\Theta\colon E\rightarrow H\F$ which is onto in homotopy. A corollary
of this is a very simple reproof of a theorem of Yagita \cite{Yag:TAMS}
that says that the map $\Theta_*\colon BP^*(BG)\rightarrow H^*(BG;\F)$
is an $F$-epimorphism. This work is described in Section~\ref{BPetc}.
Our results also show that the more interesting subalgebras $Ch_E(G)$
are going to be those arising from spectra $E$ which are either non-complex
oriented (such as $QS^0$) or else are complex oriented but periodic,
such as $K$-theory or the Johnson-Wilson theories $E(n)$.

The case of $QS^0$, already indicated, is detailed in Section~\ref{sectQS}.
In Section~\ref{exact} we examine how $Ch_E(-)$ varies as a functor when
the spectrum $E$ is changed, concentrating in particular on changing from
$E$ to an exact $E$-module spectrum, though even this simple case is more
complicated than might be thought. Sufficient is proved to show that, for
finite groups $G$, the subalgebras $Ch_{E(n)}(G)$ and $Ch_{\Ehat}(G)$ are
equal, and so the representation theoretic description of $Ch_{\Ehat}(G)$
in \cite{CCC} applies also to the subalgebras $Ch_{E(n)}(G)$. It should be
noted that the generalised character theory of \cite{HKR} applies only to
the completed spectra; this suggests that the results of \cite{CCC} ought
perhaps be independent of the work of Hopkins, Kuhn and Ravenel \cite{HKR}.

Much of the work of this paper applies to define subalgebras $Ch_E(X)$ of
$H^*(X)$ for an arbitrary space $X$, not just for spaces of the form
$X=BG$ for a finite group $G$. It is useful and in fact necessary in
some sections to work in such generality: we do so as far as possible
throughout. In the final section we concentrate on the case where $X$ has
the homotopy type of a finite CW complex. Here the spectra $E(n)$ allow
the definition of a new \lq unstable\rq\ chromatic filtration of $H^*(X)$,
compatible with the standard stable chromatic filtration of Ravenel and
others (see, for example, \cite[Def.\,7.5.3]{Rav}) defined in terms of
Bousfield localisation. For $X$ the classifying space of a finite group,
a residue of this result exists at the level of varieties.

\bigskip\noindent
{\bf Notation conventions} Throughout $p$ will denote a fixed prime and
$H^*(G)$ or $H^*(X)$ denotes the group or singular cohomology of $G$ or
$X$ respectively with constant coefficients in the field $\F$ of $p$
elements. When $E$ is a $p$-local spectrum, such as $BP$, $E(n)$ and
so on, the prime considered will be assumed to be the same one.
We identify $H^*(G)$ (group cohomology of $G$) with $H^*(BG)$ (singular
cohomology of its classifying space), generally using the former notation
except when desiring to stress the topological nature of the construction.
By the {\em variety} of an $\F$-algebra $R$, denoted $\var(R)$, we mean the
set of algebra  homomorphisms from $R$ to an algebraically closed field $k$
of characteristic $p$, topologised with the Zariski topology. In general $R$
will be graded, and as discussed in \cite{VCh} and \cite{CCC} $\var(R)$ will
be homeomorphic to $\var(R^{\rm even})$ (even when $p=2$).

\bigskip\noindent
{\bf Acknowledgements} Both authors are pleased to acknowledge the support
of the London Mathematical Society who, through a Scheme 4 grant, enabled
the second author to visit the University of Leicester where much of this
research was carried out. The first author also gratefully acknowledges
the support of a Leverhulme Research Fellowship.

\section{The basic construction}\label{basic}

We shall suppose throughout that $E$ denotes a spectrum \cite{Adams}
representing a generalised cohomology theory $E^*(-)$ and for convenience
shall suppose the coefficients $E^*=E^*($point) are concentrated in even
degrees. Moreover, we shall follow \cite{HT2,RW} and write $\E_r$ for the
$r^{\rm th}$ space in the $\Omega$ spectrum for $E$, so
$\E_r=\Omega^{\infty-r}E$ and $\E_r$ represents $E$-cohomology in degree
$r$, that is $E^r(X)=[X,\E_r]$ for any space $X$. For reasons discussed
further in \cite{CCC} it is convenient to restrict attention to the even
graded spaces $\E_{2r}$, and we shall do so throughout.
(In Section~\ref{sectQS} we restrict further to the single space $\E_0$.)
Note that each $\E_r$ is an H-space with product
$\E_r\times\E_r\rightarrow\E_r$ representing addition in $E$-cohomology
and arising from the loop space structures $\Omega\E_{r+1}=\E_r$; moreover,
each $\E_r$ is an infinite loop space.

\begin{defn}
For a spectrum $E$ and space $X$, define $Ch_E(X)$ to be the subalgebra of
$H^*(X)=H^*(X;\F)$ generated, as an $\F$-algebra, by elements of the form
$f^*(x)$ where $f\in E^{2r}(X)$ (so $f$ may be considered as a map
$f\colon X\rightarrow\E_{2r}$) and $x$ a homogeneous element in $H^*(\E_{2r})$.
\end{defn}

\begin{rem}{\em
An element of $H^k(X)$ may of course be represented by a (homotopy class of a)
map $X\rightarrow\EM_k$ where $\EM_k$ denotes an Eilenberg-Mac Lane space
of type $K(\F,k)$. It is thus equivalent to declare that $Ch_E(X)$ is the
subalgebra of $H^*(X)$ generated by maps $X\rightarrow\EM_k$ which have a
factorisation
$$X\bra{f}\E_{2r}\bra{x}\EM_k$$
for some $r$.
}\end{rem}

\begin{rem}{\em
It is immediate that this construction behaves well with respect to disjoint
unions of spaces and $Ch_E(\coprod_i X_i)=\prod_i Ch_E(X_i)$. With this
observation we shall largely restrict ourselves to connected spaces $X$.
}\end{rem}

\begin{rem}{\em
With the aid of the Cartan formula and the fact that $Ch_E(X)$ by definition
is closed under ring operations, it is not hard to check that the subalgebra
$Ch_E(X)$ is closed under the action of the Steenrod algebra.
}\end{rem}

\begin{example}\label{CCCa}{\em
(\cite{CCC} Prop.~1.6) Let $G$ be a compact Lie group and write $Ch_E(G)$ for
$Ch_E(X)$ with $X=BG$. Write $K$ for the complex $K$-theory spectrum. Then
$Ch_K(G)$ is equivalent to the classical Chern subring \cite{Thomas} of
$H^*(G)$ generated by Chern classes of irreducible representations.
}\end{example}

\begin{cor}\label{fingen}
For $G$ a finite group and $E$ a Landweber exact theory {\em
\cite{Landweber}},
$Ch_E(G)$ is a finitely generated $\F$-algebra and $H^*(G)$ is a finitely
generated $Ch_E(G)$ module.
\end{cor}

\noindent{\bf Proof}
We see from \cite[Prop.~2.8]{CCC} that if $E$ is Landweber exact and $X$ is
any
space, there is an inclusion $Ch_K(X)\subset Ch_E(X)$.
As we can identify (\ref{CCCa}) $Ch_K(G)$ with the Chern subring, over
which we know \cite{Thomas} $H^*(G)$ is finitely generated as a module, the
inclusion $Ch_K(X)\subset Ch_E(X)$ immediately tells us that $H^*(G)$ is a
finitely generated module over $Ch_E(G)$. The finite generation of $Ch_E(G)$
as an $\F$-algebra now follows from basic commutative algebra
\cite[Prop.~7.8]{AM} and the finite generation of $H^*(G)$ as
an $\F$-algebra \cite{Evens,Quillen,Venkov,Venkov2}.
\qed\bigskip

In general however, there seem at first sight to be potentially a lot of
maps $f\colon X\rightarrow\E_{2r}$ that need to be considered in the
construction of $Ch_E(X)$. The following technical lemma usefully cuts
down on othe quantity of the maps needed. For a CW complex $X$ we write
$X^{(m)}$ for its
$m$-skeleton and recall Milnor's short exact sequence
$$0\lra\lim\null\!^1 E^{2r-1}(X^{(m)})\lra E^{2r}(X)\lra\lim
E^{2r}(X^{(m)})\lra0\,.$$
Recall also that the skeleta $\{X^{(m)}\}$ define a topology (which
we shall refer to as the {\em skeletal topology\/} on the group
$\lim E^{2r}(X^{(m)})$ where the open neighbourhoods of zero are
given by
$$F_mE^{2r}(X)=\ker\left(E^{2r}(X)\lra E^{2r}(X^{(m)})\right)\,.$$

\begin{lem}\label{topgen}
Suppose $X$ is a CW complex. Then $Ch_E(X)$ is generated by elements
$f^*(x)$ where $f\colon X\rightarrow\E_{2r}$ run over a set of topological
generators of the $E^*$-module $E^{\rm even}(X)$ modulo phantom maps.
\end{lem}

\noindent{\bf Proof}
First note that we can disregard any phantom map $f\colon
X\rightarrow\E_{2r}$.
Suppose such a map gave rise to an element $f^*(x)\in H^k(X)$ via some element
$x\in H^k(\E_{2r})$. As $H^k(X)\cong H^k(X^{(m)})$ for any $m>k+1$ with the
isomorphism induced by restriction $\rho\colon X^{(m)}\rightarrow X$, the
element $f^*(x)$ is non-zero if and only if the composite
$$X^{(m)}\bra{f\rho}\E_{2r}\bra x\EM_k$$
is non-zero.
However, by definition of $f$ being a phantom map, $f\rho$ is trivial and
we are just picking up
the zero element. We can thus restrict the elements $f$ we need to those
having non-zero image in $\lim E^{\rm even}(X^{(m)})$.

Next observe that is is sufficient to consider maps
$f\colon X\rightarrow\E_{2r}$ which are $E^*$-module generators of
$E^{\rm even}(X)$. The $E^*$-action on $E^*(X)$ may be represented by maps
$$\mu\colon E^{2t}\times\E_{2r}\rightarrow\E_{2(r+t)}$$
where the group $E^{2t}$ is regarded as a space with the discrete topology.
Now consider an element of $Ch_E(X)$ given by $(\alpha f)^*(x)$ where
$f\in E^{2r}(X)$ and $\alpha\in E^{2t}$. The map $\alpha f$ is represented
by the composite
$$X\bra\Delta X\times X\bra{\alpha\times f}E^{2t}\times
\E_{2r}\bra\mu\E_{2(r+t)}$$
but this is equivalent to a composite
$$X\bra g\{\alpha\}\times\E_{2r}\cong\E_{2r}\bra\mu\E_{2(r+t)}$$
where we write $\mu$ also for the restriction of $\mu$ to this specific
component. Thus we have the equality
$$(\alpha f)^*(x)=g^*(x\mu)\,.$$

\noindent
Meanwhile, addition in $E^{2r}(X)$ is represented by a map
$$\sigma\colon\E_{2r}\times\E_{2r}\lra\E_{2r}\,.$$
Consider an element of $Ch_E(X)$ given by $(f_1+f_2)^*(x)$ with
$f_i\in\E_{2r}$ and $x\in H^t(X)$. Then this element can be written as
the composite
$$X\bra\Delta X\times X\bra{f_1\times
f_2}\E_{2r}\times\E_{2r}\bra\sigma\E_{2r}\bra x\EM_t\,.$$
If we write $\sigma^*(x)$ as $\sum x'\otimes x''\in H^*(\E_{2r})\otimes
H^*(\E_{2r})$ this shows that $(f_1+f_2)^*(x)=\sum f_1^*(x')f_2^*(x'')$,
an element in the ring generated using just $f_1$ and $f_2$.

Finally, we show that we need only consider topological generators of
$E^{\rm even}(X)$ mod phantom maps. Once again consider a generator
$f^*(x)\in Ch_E(X)$ and let $x\in H^k(\E_{2r})$. Suppose $f\in E^{2r}(X)$
is represented by a sequence of elements $f_m\in E^{2r}(X)$, $m=1,2,\ldots$,
where $f-f_m\in F_m E^{2r}(X)$ and $f_m$ is in the (algebraic) $E^*$-span of
a given set of topological generators. As $H^k(X)\cong H^k(X^{(m)})$ for all
$m>k+1$, the elements $f^*(x)$ and $(f_m)^*(x)$ are identical for all
$m>k+1$. This concludes the proof.
\qed

\begin{rem}{\em
If $E$ is a ring spectrum and so $E^*(X)$ is an $E^*$-algebra, a similar
argument to the one above for sums shows that we can reduce the maps
$f\in E^*(X)$ further to include only topological ring generators of
$E^*(X)$ mod phantoms.
}\end{rem}

\begin{warning}{\em
The result (\ref{topgen}) does not imply that the $Ch_E(X)$ may be computed
for a CW complex $X$ in terms of $Ch_E(X^{(m)})$. While an element of
$Ch_E(X)$ certainly gives rise, by restriction, to an element of
$\lim Ch_E(X^{(m)})$, it is possible to have a tower of elements in
$\lim Ch_E(X^{(m)})$ which do not lift to $Ch_E(X)$.
}\end{warning}

We conclude this section with some remarks about the relation of the
subalgebra $Ch_E(X)$ to the Bousfield $E$-localisation of $X$.  Recall
\cite{bousfield} defines for a generalised homology theory $E_*(-)$ a
localisation functor $L_E$ with a natural transformation
$\eta_X\colon X\mapsto L_EX$ satisfying certain universal properties.

\begin{prop}\label{Bloc}
For any space $X$, the subalgebra $Ch_E(X)\subset H^*(X)$ is contained
in the image of $\eta_X^*\colon H^*(L_EX)\rightarrow H^*(X)$.
\end{prop}

\noindent{\bf Proof}
The spaces $\E_{2r}$ are themselves $E$-local. Hence any map
$f\colon X\rightarrow \E_{2r}$ factors through the localisation
$\eta_X\colon X\mapsto L_EX$. Thus any generator of $Ch_E(X)$ lies
in $\eta_X^*$. As $\eta_X^*$ is an algebra homomorphism it is closed
under sums and products and the result follows.
\qed

\begin{rem}{\em
The inclusion $Ch_E(X)\subset\Image(\eta_X^*)$ of (\ref{Bloc}) will
usually be strict. For example, if $E$ is any Landweber exact spectrum,
$H^*(\E_{2r})$ lies entirely in even dimensions \cite{BH} and hence
$Ch_E(X)\subset H^{\rm even}(X)$; there are many examples where
$\Image(\eta_X^*)$ contains odd dimensional elements, for example if we
take $X=B\Z/p$ and $E=BP$.

However, even the inclusion $Ch_E(X)\subset\Image(\eta_X^{\rm even})$ for
such spectra will also likely be strict. For an example here, let $X$ be
the classifying space of an elementary abelian group of rank $3$ then
(as $X$ is then a nilpotent space) \cite{bousfield} tells us that $L_{BP}X$
is just the $p$-localisation of $X$ and so $\Image(\eta_X^{\rm even})$ is
the whole of $H^{\rm even}(X)$. We shall see by Corollary~\ref{stabunstab}
that $Ch_{BP}(X)$ can be identified with the image of $BP^*(X)$ under the
Thom map, and for this particular $X$ this image is not the whole of
$H^{\rm even}(X)$.
}\end{rem}

\section{Stabilisation and $BP$-theory}\label{BPetc}

There is no reason to suppose that the maps $x\colon\E_{2r}\rightarrow\EM_k$
used in the definition of $Ch_E(X)$ commute with the loop space structure on
these spaces. This observation allows us to make the following definition.

\begin{defn}
Define the subalgebra $\Sigma^d Ch_E(X)$ of $Ch_E(X)\subset H^*(X)$ to be
that generated by elements $f^*(x)$ where $x\colon \E_{2r}\rightarrow\EM_k$
is a $d$-fold loop map, {\em i.e.}, there is map
$y\colon\E_{2r+d}=\Omega^{-d}\E_{2r}\rightarrow\EM_{k+d}=\Omega^{-d}\EM_k$ with
$x=\Omega^dy$. Likewise, define $Ch_E^s(X)$, the stable
$E$-Chern subalgebra, to be the intersection of all $\Sigma^d Ch_E(X)$,
that is the subalgebra generated by elements $f^*(x)$ where
$x\colon \E_{2r}\rightarrow\EM_k$ is a map of infinite loop spaces.
\end{defn}

\begin{rem}\label{KStriv}{\em
The process of stabilisation may yield only the trivial subalgebra. For
example, the $K$-theory subalgebra $Ch_K(X)$ is in general non-trivial (for
example, taking $X$ to be $B\Z/p$, we have $Ch_K(B\Z/p)=H^{\rm even}(B\Z/p)$).
However, there are no non-trivial stable maps $K\rightarrow H\F$ and thus
$Ch_K^s(X)=0$ for all $X$. Similar remarks apply with $K$ replaced by any
of the $E(n)$ spectra.
}\end{rem}

Now suppose the spectrum $E$ is a ring spectrum and has a Thom map, by
which for the moment we just mean a distinguished map of ring spectra
$\Theta\colon E\rightarrow H\F$ which is onto in homotopy. Our principal
examples are $E=BP$ or $MU$, but we also include spectra such as $k(n)$,
$P(n)$ and $BP\langle n\rangle$; the periodic spectra such as $K(n)$
and $E(n)$ are obviously not included. We write $\Theta$ also for the
corresponding individual infinite loop maps $\E_{2r}\rightarrow\EM_{2r}$.

\begin{defn}\label{thetadef}
Define the subalgebra $Ch_E^\Theta(X)$ of $H^*(X)$ to be that generated by
elements of the form $f^*(\Theta)$ where $f\in E^{2r}(X)$ and $\Theta$ is the
appropriate infinite loop map, as
just given. Of course this is just the image of
$E^*(X)$ in $H^*(X)$ under the ring homomorphism $\Theta_*$. In this
situation we have natural inclusions
$$Ch_E^\Theta(X)\subset Ch_E^s(X)\subset Ch_E(X)\subset H^*(X).$$
\end{defn}

Before the next result let us remark again that until otherwised mentioned,
all spectra such as the $E$ and $T$ following are supposed to have
coefficients
 concentrated in even degrees. We recall the notion of Landweber exact spectra
\cite{Landweber}. This class of spectra include complex cobordism $MU$ and
the Brown-Peterson theories $BP$, the Johnson-Wilson spectra $E(n)$ and their
various completions such as $\Ehat$ and Morava $E$-theory, it also contains
complex $K$-theory and many examples of elliptic spectra. The topology of the
spaces in the $\Omega$ spectrum for a Landweber exact theory is well
understood
\cite{BH, HH, hun, HT2, RW}; in particular, $H_*(\E_{2r};\Z)$ is torsion free
and concentrated in even dimensions.

\begin{prop}\label{lift}
Let $E$ be any Landweber exact spectrum and suppose $T$ to be a ring spectrum
with Thom map. Then for any space $X$ there is an inclusion of subalgebras
$$Ch_E(X)\subset Ch^\Theta_T(X)\,.$$
\end{prop}

\noindent{\bf Proof}
The Thom map $\Theta\colon T\rightarrow H\F$ induces a map of
Atiyah-Hirzebruch spectral sequences for each space $\E_{2r}$
which on $E_2$-pages is a surjection
$$H^*(\E_{2r};T^*)\lra H^*(\E_{2r};\F)\,.$$
By the properties of $H_*(\E_{2r};\Z)$ for a Landweber exact spectrum $E$
mentioned above, both these spectral sequences collapse and we conclude that
any map $x\colon\E_{2r}\rightarrow\EM_k$ can be lifted through the Thom map
$\E_{2r}\bra l\underline{T\!}\,_k\bra\Theta\EM_k$. (Note that only even
values of $k$ arise as $H^*(\E_{2r})$ is concentrated in even dimensions.)
Thus any element $f^*(x)\in Ch_E(X)$ with $f\colon X\rightarrow\E_{2r}$ and
$x\in H^k(\E_{2r})$ may be written as $(lf)^*(\Theta)\in Ch_T^\Theta(X)$.
\qed\bigskip

The spectra $MU$ and $BP$ each satisfy the hypotheses on both $E$ and $T$
in Proposition~\ref{lift} (necessarily assuming the prime $p$ of the
underlying field $\F$ is the same as that for the version of $BP$ considered).
As $Ch_T^\Theta(X)\subset Ch_T(X)$, by choosing both $E$ and $T$ to be $BP$
(and likewise choosing both to be $MU$), (\ref{lift}) gives us

\begin{cor}\label{stabunstab}
For any space $X$ we have equalities
\begin{alignat}{2}
Ch_{BP}^\Theta(X)&=&Ch_{BP}^s(X)&=Ch_{BP}(X) \notag\\
Ch_{MU}^\Theta(X)&=&Ch_{MU}^s(X)&=Ch_{MU}(X)\,.\tag*{$\square$}
\end{alignat}
\end{cor}

\begin{rem}{\em
This result should be contrasted with Remark~\ref{KStriv} which, with
Corollary~\ref{fingen} observes that for $E$ any of the periodic
Landweber exact spectra listed above, the subalgebras $Ch_E(X)$ and
$Ch_E^s(X)$ are quite different -- the former, for $X=BG$, containing
the whole Chern subring, the latter being trivial.
}\end{rem}

\noindent
Choosing one of $E$ or $T$ to be $BP$ and the other to be $MU$,
Proposition~\ref{lift} tells us

\begin{cor}
For $X$ any space, $Ch_{BP}(X)=Ch_{MU}(X)$.\qed
\end{cor}

\begin{thm}\label{ThomThm}
Let $T$ be a complex oriented ring spectrum with a Thom map and suppose
$G$ to be a finite group. Then the inclusions $Ch_T(G)\subset H^*(G)$
and $Ch_T^\Theta(G)\subset H^*(G)$ are $F$-isomorphisms. Hence, by the
first observation in
{\em(\ref{thetadef})}, for such $T$ the map
$\Theta_*\colon E^*(BG)\rightarrow H^*(G)$ is an $F$-epimorphism.
\end{thm}

\noindent
We offer two proofs of this result, both resting on the stabilisation results
just proved. The first proof uses also the main theorem of \cite{CCC}, while
the second uses instead work of Carlson \cite{carlson} or of Green and Leary
\cite{VCh} on corestriction of Chern classes.

\bigskip\noindent
{\bf First Proof} By Proposition~\ref{lift}, for every $n=1,2,\ldots$,
there are inclusions
$$Ch_{\Ehat}(G)\subset Ch_T^\Theta(G)\subset Ch_T(G)\subset H^*(G)\,.$$
However, by \cite[Theorem 0.2]{CCC} the inclusion $Ch_{\Ehat}(G)\subset H^*(G)$
is an $F$-iso\-mor\-phism if $n$ is not less than the $p$-rank of $G$.
\qed

\bigskip\noindent
{\bf Second Proof} By Proposition~\ref{lift} there are inclusions
$$Ch_K(G)\subset Ch_T^\Theta(G)\subset Ch_T(G)\subset H^*(G)\,.$$
We follow \cite{VCh} and write $\overline{R(G)}$ for the closure of a
functorial subring $R(G)\subset H^*(G)$ under corestriction of elements
from $R(H)$ for all subgroups $H$ of $G$. By \cite[\S8]{VCh} or \cite{carlson}
(see \cite{VCh} for detailed discussion), $\overline{Ch_K(G)}$ is
$F$-isomorphic to $H^*(G)$, implying $\overline{Ch_T^\Theta(G)}$ is
also $F$-isomorphic to $H^*(G)$. However, as corestriction is a stable
construction, $\overline{Ch_T^\Theta(G)}=Ch_T^\Theta(G)$ and so
$Ch_T^\Theta(G)$ and hence $Ch_T(G)$ is $F$-isomorphic to $H^*(G)$.
\qed

\begin{rem}{\em
If we take $T=BP$ or $k(n)$ in Theorem~\ref{ThomThm} we recover Yagita's
results \cite[(4.2), (4.3)]{Yag:TAMS}. As noted at the end of the previous
section, the inclusions $Ch_{BP}(G)\subset H^{\rm even}(G)$ can however be
strict, for example when $G$ is elementary abelian of rank $3$. In the next
section we shall show that the result (\ref{ThomThm}) can fail when $T$ is
not complex oriented.
}\end{rem}

\begin{cor}
Suppose $T$ is a complex oriented ring spectrum with a Thom map and $G$ is
a finite group. Let $\mathcal{A}(G)$ be the category of elementary abelian
subgroups of $G$ with morphisms those inclusions $V\rightarrow W$ given by
conjugation by an element of $G$. Then there is a homeomorphism of varieties
$$\var(Ch_T(G))\cong\colim_{V\in{\mathcal{A}}(G)}\var(H^*(V))\,.$$
\end{cor}

\noindent{\bf Proof}
This follows from Theorem~\ref{ThomThm} which identifies $\var(Ch_T(G))$
with $\var(H^*(G))$ and Quillen's theorem \cite{Quillen} which describes
$\var(H^*(G))$ as the colimit of the $H^*(V)$ over the category indicated.
\qed

\section{Subalgebras defined by $QS^0$}\label{sectQS}

Suppose throughout this section that $G$ is a finite group. The main result
below identifies the subalgebra $\CS(G)$ of $H^*(G)$ by which we mean the
subalgebra of $H^*(G)$ generated by the images in cohomology of maps
$BG\rightarrow QS^0$.
Recall that $QS^0=\dlim_{n\rightarrow\infty}\Omega^n\Sigma^n S^0$ and is the
infinite loop space that represents stable cohomotopy in degree zero,
{\em i.e.}, $\pi^0_s(X)=[X,QS^0]$. We restrict attention here to the zero
space $QS^0$ as, by the Segal conjecture \cite{carlsson}, there are no
non-trivial maps $BG\rightarrow QS^n$ for $n>0$. The cases $n<0$ present
formidable difficulties.

We begin by recalling the algebra $S_h=S_h(G)\subset H^*(G)$ of
\cite[Def.\ 1.2]{GLS}. For a finite $G$-set $X$ of cardinality $n$,
a choice of bijection between $X$ and the set $\{1,\ldots,n\}$ induces
a homomorphism $\rho_X\colon G\rightarrow \Sigma_n$ where $\Sigma_n$ is
the symmetric group on $n$ letters. It thus induces an algebra homomorphism
$\rho_X^*\colon
H^*(\Sigma_n)\rightarrow H^*(G)$. For a fixed space $X$ two choices of
$\rho_X$ differ only by an inner automorphism of $\Sigma_n$ and
so the homomorphism $\rho_X^*\colon H^*(\Sigma_n)\rightarrow H^*(G)$
depends only on $X$ and not on the choice of bijection.
Then $S_h=S_h(G)$ is defined as the subalgebra of $H^*(G)$ generated by
elements of
the form $\rho_X^*(x)$ as $X$ runs over all finite $G$ sets and where the
$x$ are homogeneous
elements of $H^*(\Sigma_n)$.

\begin{thm}\label{QS0thm}
For $G$ finite there is an algebra isomorphism $\CS(G)\cong S_h(G)$.
\end{thm}

\noindent{\bf Proof}
First we show that $S_h(G)$ is a subalgebra of $\CS(G)$. By definition,
$S_h(G)$ is generated by elements of the form $f^*(x)$ where $f$ is drawn
from some particular class of maps $BG\rightarrow B\Sigma_n$ and
$x\in H^*(\Sigma_n)$. Given such an element, consider the composite
$$BG\bra f B\Sigma_n\bra{i_n}B\Sigma_\infty\bra D QS^0$$
where $i_n$ is induced by the incusion of $\Sigma_n$ in the infinite
symmetric group, and $D$ is the Dyer-Lashof map. By \cite{Priddy} $D$
induces an isomorphism in cohomology, while by Nakaoka \cite[\S7]{nakaoka}
the induced map $i_n^*\colon H^*(B\Sigma_\infty)\rightarrow H^*(B\Sigma_n)$
is a surjection.
Thus there is an element $y\in H^*(QS^0)$ with $x=(Di_n)^*(y)$. Hence
$f^*(x)=(Di_nf)^*(y)\in\CS(G)$.

To show the converse we use the Segal conjecture \cite{carlsson}. Recall
that this identifies the homotopy classes of maps $BG\rightarrow QS^0$ with
the elements of the Burnside ring, completed at the augmentation ideal.
Specifically, to an element of the Burnside ring $A(G)$ we can associate a
finite $G$-set $X$; the Segal conjecture tells us that the map which sends
this element to the resulting composite
$$BG\bra f B\Sigma_n\bra{i_n}B\Sigma_\infty\bra D QS^0$$
(where again, $n$ is the cardinality of $X$) is continuous (with respect to
the topologies induced by powers of the augmentation ideal in $A(G)$ and the
skeleta in $\pi_s^0(BG)$) and an isomorphism after completion; the ring
$\pi_s^0(BG)$ is already complete under the skeletal topology. Thus we can
take {\em topological\/} generators of $\pi_s^0(BG)$ from the ring $A(G)$
itself and any element $g^*(y)\in\CS(G)$ represented by such an element
$g\in\pi_s^0(BG)$ and $y\in H^*(QS^0)$ can be realised as an element
$f^*(x)\in S_h(G)$ where $x=(Di_n)^*(y)$. Lemma~\ref{topgen} now shows that
this is sufficient to prove that $\CS(G)$ is contained in $S_h(G)$.
\qed

\begin{rem}{\em
A representation theoretic description of the variety associated to $S_h(G)$
is given in \cite{GLS}. The work there combined with Theorem~\ref{QS0thm}
shows there is a natural homeomorphism
$$\var(\CS(G))\cong\colim_{V\in{\mathcal{A}}_h(G)}\var(H^*(V))$$
where $\mathcal{A}_h(G)$ denotes the category whose objects are elementary
abelian subgroups of $G$ and whose morphisms are those injective group
homomorphisms $f\colon V\rightarrow W$ for which $f(U)$ is conjugate in $G$
to $U$ for every subgroup $U$ of $V$.
}\end{rem}

\begin{rem}{\em
It is also noted in \cite{GLS} that in general neither $S_h(G)$ nor $Ch(G)$,
the classical Chern subring, are contained in the other. See \cite{GLS} for
worked examples.
}\end{rem}

\section{Subalgebras defined by exact spectra}\label{exact}

In this section we examine the relation between $Ch_E(X)$ and $Ch_F(X)$ for
spectra $E$ and $F$ that are themselves related; specifically,
we are interested in the case where $F$ is {\em exact\/} over $E$, by which
we mean that $E$ is a ring spectrum, $F$ an $E$-module spectrum and the
corresponding homology theories are related by a natural equivalence
$$F_*(X)\cong F_*\tensor_{E_*}E_*(X)$$
for any space $X$.

It is of course not true that if $F$ is exact over $E$ then $Ch_E(X)$ and
$Ch_F(X)$ are identical: for example, complex $K$-theory is exact over
$BP$-theory, but the inclusion $Ch_K(X)\subset Ch_{BP}(X)$ is certainly
strict in many cases. However, exactness of $F$ over $E$ is sufficient,
under suitable finiteness conditions, to give in general an inclusion
$Ch_F(X)\subset Ch_E(X)$; see the results (\ref{fiveB}), (\ref{fiveC}),
(\ref{fiveD}) and (\ref{fiveH}) below.

Our main goal here though, Theorem~\ref{thm:five}, is to show that the
complete theories $\Ehat$ used in \cite{CCC} are not the only ones that
can give the topological interpretation of the subgroup categories
$\mathcal{A}^{(n)}$ of \cite{VCh,CCC}, and that the more classical
Johnson-Wilson theories $E(n)$ \cite{JW} satisfy a similar theorem.
The Baker-W\"urgler completions $\Ehat$, or equivalently the Morava
$E$-theories, were necessary in \cite{CCC} so as to apply the generalised
character theory of Hopkins, Kuhn and Ravenel \cite{HKR} and the
completeness assumptions seem integral
to that character theory. Thus Theorem~\ref{thm:five} below has apparently
by-passed the subtleties of \cite{HKR} and suggests it ought to be possible
to find a proof of Theorem~\ref{thm:five} without reference to \cite{HKR};
it would be interesting to see such a direct argument.

Throughout this section we shall not only continue the convention that any
subalgebra $Ch_E(X)$ or $Ch_F(X)$ is defined using just the even graded spaces
in the $E$ or $F$ $\Omega$-spectra, but we shall further restrict these
subalgebras to those lying in even
dimensional cohomology, that is, by $Ch_E(X)$ we shall really mean
$Ch_E(X)\cap H^{\rm even}(X)$,  {\em
etc}. This convention is made for notational simplicity and could be
avoided if one wanted. However, as the
main examples we wish to discuss here are of Landweber exact spectra,  we
see by \cite{BH} or by
\cite{HT2,RW} that this is not actually a restriction at all: there are no
odd dimensional elements in the
original subalgebras
$Ch_E(X)$ for such $E$.

We begin with the following general observation about the subalgebras
$Ch_E(X)$ and $Ch_F(X)$ for arbitrary spectra $E$ and $F$ (not necessarily
exact or even one a module spectrum over the other) arrived at by considering
the universal example of $X=\FF_{2*}$.

\begin{prop}\label{fiveA}
Suppose the cohomology of each space $H^{\rm even}(\FF_{2r})$ is of
finite type and suppose $Ch_E(\FF_{2r})=H^{\rm even}(\FF_{2r})$ for each
$r$. Then $Ch_F(X)\subset Ch_E(X)$ for every space $X$.
\end{prop}

\noindent{\bf Proof}
Consider a generator $f^*(x)$ of $Ch_F(X)$, where $f\colon
X\rightarrow\FF_{2r}$ and $x\in H^{2t}(\FF_{2r})$. As
$Ch_E(\FF_{2r})=H^{\rm even}(\FF_{2r})$, the
finiteness assumption means that we can write $x$ as a finite sum of
products of elements $g_i^*(y_i)\in Ch_E(\FF_{2r})$,
with $g_i\colon\FF_{2r}\rightarrow\E_{2s_i}$ and $y_i\in H^*(\E_{2s_i})$.
The element $f^*(x)$ may thus be written as the composite
$$X\bra f\FF_{2r}\bra\Delta\prod_i\FF_{2r}\bra{\prod g_i}
\prod_i\E_{2s_i}\lra\prod_i\EM_{2t_i}\lra\EM_{2t}$$
where the last arrow denotes the summing and multiplication operations.
The composite of the final two arrows denotes an element of
$H^{2t}(\prod_i\E_{2s_i})$ and as such we can write, using the
K\"unneth theorem for mod $p$ cohomology, in the form
$\sum y^{(1)}\otimes\cdots\otimes y^{(q)}$ (where the sum is potentially
infinite). If we write $f_i$ for
the projection of the composite $X\bra f\FF_{2r}\lra\prod_i\E_{2s_i}$
to the $i^{\rm th}$ factor, we see that
$$f^*(x)=\sum f^*_1(y^{(1)})\cdots f^*_q(y^{(q)})\, ,$$
an element of $Ch_E(X)$: note that this sum will contain only a finite
number of non-zero terms by virtue of
the finite type assumption on $H^{\rm even}(\FF_{2r})$, even if the first
sum $\sum
y^{(1)}\otimes\cdots\otimes y^{(q)}\in H^*(\prod_i\E_{2s_i})$ is infinite.
\qed

\begin{prop}\label{fiveB}
Suppose $F$ is an exact $E$-module spectrum and that the space
$X$ is of finite type. Then $Ch_F(X)\subset Ch_E(X)$.
\end{prop}

\noindent{\bf Proof}
Assume $X$ is connected, for else we can argue by connected components.
If $X$ had the homotopy type of a finite CW complex then
$F^*(X)=F^*\tensor_{E^*}E^*(X)$; if $X$ is only of finite type then this
statement is only true after completion of the tensor product. The finite
type hypothesis however tells us that the completion is with respect to the
skeletal topology (as in Section~\ref{basic}):
$$F^*(X)=F^*\widehat{\tensor_{E^*}}E^*(X)=
\lim_{m\rightarrow\infty}\left(F^*\tensor_{E^*}E^*(X^{(m)})\right)\,.$$
Thus topological generators of $F^*(X)$ may be taken from elements of the
(uncompleted) tensor product $F^*\tensor_{E^*}E^*(X)$. Such an element,
considered as a map $X\rightarrow\FF_{2r}$, may thus be lifted as a finite sum
of maps through spaces $F^*\times\E_{2s_i}$, where $F^*$ is regarded as a
discrete space, and thus we can represent $f$ as a composite
$$X\lra\prod_i\left(F^*\times\E_{2s_i}\right)\lra\FF_{2r}\,.$$
By the connectedness assumption, such a map actually lifts through a finite
product of spaces $\E_{2s_i}$.

Now consider a generator $f^*(x)\in Ch_F(X)$ as usual. Factoring $f$ as
above as a map $X\lra\prod_i\E_{2s_i}\lra\FF_{2r}$ together with the finite
type assumption on $X$ now shows
$f^*(x)$ to be in the subalgebra $Ch_E(X)$.
\qed

\begin{cor}\label{fiveC}
If $F$ is an exact $E$-module spectrum and the spaces $\FF_{2r}$ representing
$F$ are of finite type, then $Ch_F(X)\subset Ch_E(X)$ for every space $X$.
\end{cor}

\noindent{\bf Proof} 
We apply Proposition~\ref{fiveB} with $X=\FF_{2r}$ and deduce that
$$Ch_F(\FF_{2r})\subset Ch_E(\FF_{2r})\subset H^*(\FF_{2r}).$$ 
But
$Ch_F(\FF_{2r})=H^*(\FF_{2r})$ by definition of $Ch_F(-)$ and hence
$Ch_E(\FF_{2r})=H^*(\FF_{2r})$. The result now follows from
Proposition~\ref{fiveA}.
\qed

\begin{rem}\label{fiveD}{\em
There are a number of examples to which this corollary may be applied.
First of all, we could take $E=MU$ and $F=K$ and we recover the inclusion
$Ch_K(X)\subset Ch_{MU}(-)$ of Propositions~\ref{CCCa} and~\ref{lift}. New
results however include the inclusions
\begin{alignat*}{3}
&Ch_{KO}(X)\subset Ch_{MSp}(X)&\quad
&\!Ch_{KO}(X)\subset Ch_{MSpin}(X)&\quad
&\!Ch_{KO}(X)\subset Ch_{El}(X)\\
&Ch_{K}(X)\subset Ch_{MSpin^c}(X)&\quad
&Ch_{El}(X)\subset Ch_{MSpin}(X)&
\end{alignat*}
where $MSpin^c$ is the self-conjugate spin cobordism of \cite{HopHov} and
$El$ is the (non-complex oriented) {\em integral\/} elliptic theory of Kreck
and Stolz \cite{krst:hpt}; see \cite{Hov}. With care, the arguments above
may be extended to prove also the inclusion
$$Ch_{K(n)}(X)\subset Ch_{P(n)}(X)$$
using Yagita's $I_n$ version of the Landweber exact functor
theorem \cite{Yag:Land}.
}\end{rem}

\begin{lem}\label{fiveE}
Suppose $E$ and $F$ are spectra equipped with a map of spectra $\phi\colon
E\rightarrow F$
inducing an inclusion in mod $p$ homology of the corresponding
$\Omega$-spectra,
$H_*(\E_*)\hookrightarrow H_*(\FF_*)$. Then $Ch_E(X)\subset Ch_F(X)$ for
any space $X$.
\end{lem}

\noindent{\bf Proof}
If $\phi$ induces an injection in homology, it induces a
surjection $H^*(\FF_*)\onto H^*(\E_*)$ in cohomology. Thus, given
a generator $f^*(x)$ of $Ch_E(X)$ with $f\colon X\rightarrow\E_{2r}$
and $x\in H^{2t}(\E_{2r})$,  the element $x$ may be lifted as $x=\phi^*(y)$
for some $y\in H^{2t}(\FF_{2r})$ and $f^*(x)$ can be represented as a composite
$$X\bra f\E_{2r}\bra\phi\FF_{2r}\lra\EM_{2t}\,,$$
an element of $Ch_F(X)$.
\qed

\bigskip\noindent
This result will give us a converse inclusion to that offered by the previous
results. Of course there is no reason to suppose that even if $F$ is exact
over
$E$ there should be an inclusion of the homology of theory $\Omega$-spectra.
For example, $K$ is an exact $BP$-module spectrum but the map
$H_*(\underline{BP\!}\,_*)\rightarrow H_*(\underline{K\!}\,_*)$
is not an inclusion. Thus we cannot conclude from (\ref{fiveD}) the
equality of $Ch_{BP}(-)$ and $Ch_K(-)$, but then, as noted above, they are
not equivalent.

We can however apply these results to the case of $E=E(n)$ and $F=\Ehat$: we
know \cite{bwa,bwb} that $\Ehat$ is an exact spectrum over $E(n)$ and it is
shown in \cite{hun} that the map $H_*(\underline{E(n)\!}\,_*)\rightarrow
H_*(\underline{\Ehat\!}\,_*)$ induced by the completion is an inclusion. We
thus obtain
by (\ref{fiveB}) and (\ref{fiveE})

\begin{cor}\label{fiveF}
For $X$ with the homotopy type of a finite type CW complex,
there is an equivalence
$$Ch_{E(n)}(X)=Ch_{\Ehat}(X)\,.$$
\qed
\end{cor}

\noindent
Recall the categories $\mathcal{A}^{(n)}=\mathcal{A}^{(n)}(G)$ of elementary
abelian subgroups of a group $G$ defined in \cite{VCh}: the objects of
$\mathcal{A}^{(n)}(G)$ are the
elementary abelian subgroups of $G$, denoted $W$, $V$, {\em etc.}, and the
morphisms from $W$ to $V$ are
those injective group homomorphisms $f\colon W\rightarrow V$ with the
property that if $U$ is any subgroup
of $W$ of rank at most $n$, the restriction of $f$ to $U$ is given by
conjugation by some element of $G$.
Corollary~\ref{fiveF}, together with the main result of
\cite{CCC}, now gives

\begin{thm}\label{thm:five}
For a finite group $G$ there is a homeomorphism of varieties
$$\Vark{Ch_{E(n)}(G)}\cong\colim_{V\in\An(G)} \Vark{H^*(V)}\,.$$\qed
\end{thm}

\noindent
We can now also explain the relationship between $Ch_{E(1)}(X)$ and $Ch_K(X)$.
For $X$ the classifying space of a finite group it was shown in \cite{CCC}
that $Ch_{\widehat{E(1)}}(X)$ and $Ch_K(X)$ are $F$-isomorphic, the latter
being the classical Chern subring. The former, by (\ref{fiveF}), is
equivalent to $Ch_{E(1)}(X)$.

\begin{prop}\label{fiveH}
For any space $X$ the subalgebras $Ch_{E(1)}(X)$ and $Ch_K(X)$ of $H^*(X)$
are equal.
\end{prop}

\noindent{\bf Proof}
Write $K\Z_{(p)}$ for $p$-local complex $K$-theory; the localisation
map $\Z\rightarrow\Z_{(p)}$ induces a map of ring spectra
$K\bra l K\Z_{(p)}$ and in fact $K\Z_{(p)}$ is exact over $K$ via this map.
Moreover, a routine calculation shows that the induced map
$H_*(\underline{K\!}\,_*)\bra{l_*}H_*(\underline{K\Z_{(p)}\!}\,_*)$ is an
inclusion (for example, the calculations of the form of \cite[\S3]{HR} show
that in positive dimensions this is an isomorphism). Thus it follows from
Corollary~\ref{fiveC} and Lemma~\ref{fiveE} that $Ch_K(X)=Ch_{K\Z_{(p)}}(X)$
for all $X$. However, $K\Z_{(p)}$ splits into a sum of suspensions of copies
of $E(1)$ and in fact is exact over $E(1)$; moreover, the splitting gives the
necessary inclusion of homology to allow us to use (\ref{fiveC}) and
(\ref{fiveE}) again to prove that $Ch_{E(1)}(X)=Ch_{K\Z_{(p)}}(X)$.
\qed

\section{Unstable chromatic filtrations}\label{finite}

In \cite{CCC} it was shown for a finite group $G$ that the set of subalgebras
$Ch_{\Ehat}(G)$ (which by Corollary~\ref{fiveF} we can now identify with
$Ch_{E(n)}(G)$) formed an \lq\lq$F$-filtration\rq\rq\ of $H^*(G)$, in the
sense that their associated varieties formed a sequence of quotient spaces
$$\Vark{H^*(G)}\onto\cdots\onto \Vark{\Ch_{\widehat{E(n+1)}}(G)}\onto
\Vark{\Ch_{\widehat{E(n)}}(G)}\onto\cdots\onto k\,.$$
The exact relation between the subalgebras $Ch_{E(n)}(G)$ and
$Ch_{E(n+1)}(G)$ is however uncertain: there is no reason to suppose
that there is an actual inclusion $Ch_{E(n)}(G)\subset Ch_{E(n+1)}(G)$.

In this section we examine the subalgebras $Ch_{E(n)}(X)$ for a space $X$
with the homotopy type of a finite CW complex. Here we have a stronger
result which shows the subalgebras $Ch_{E(n)}(X)$ form an actual filtration
of the algebra $H^{\rm even}(X)$; by the work of Section~\ref{basic} this is
as subalgebras over the Steenrod algebra and compatible with the (unstable)
chromatic filtration of the space $X$ ({\em cf.} \cite[Def.\,7.5.3]{Rav}
and see Remark~\ref{six:rem} below).
Our main result is as follows.

\begin{thm}\label{six}
Let $X$ be a finite CW complex. Then for every $n=1,2,\ldots$ there is an
inclusion $Ch_{E(n)}(X)\subset Ch_{E(n+1)}(X)$ as subalgebras of
$H^{\rm even}(X)$.
\end{thm}

\noindent{\bf Proof}
We use the associated spectra $BP\langle n\rangle$ and Wilson's
Splitting Theorem of \cite{wilson2}.
Specifically, the spectrum $BP\langle n\rangle$ has homotopy
$\Z_{(p)}[v_1,\ldots,v_n]$ with $v_i$ of dimension $|v_i|=2(p^i-1)$,
and for $r<(p^n+\cdots+p+1)$
$$\underline{BP\langle n+1\rangle\!}\,_{2r}\simeq
\underline{BP\langle n\rangle\!}\,_{2r}
\times\underline{BP\langle n+1\rangle\!}\,_{2r+|v_{n+1}|}$$
as H-spaces \cite[Corollary 5.5]{wilson2}. Recall also that we can define
the spectrum $E(n)$ as $v_n^{-1}BP\langle n\rangle$ and the space
$\underline{E(n)\!}\,_{2r}$ can be constructed as the homotopy colimit of
$$\underline{BP\langle n\rangle\!}\,_{2r}\bra{v_n}
\underline{BP\langle n\rangle\!}\,_{2r-|v_n|}\bra{v_n}
\underline{BP\langle n\rangle\!}\,_{2r-2|v_n|}\lra\cdots$$
where the maps, as indicated, represent in homotopy multiplication by $v_n$.

Now consider a generator of $Ch_{E(n)}(X)$ which we take to be represented
by $f^*(x)$ for a map $f\colon X\rightarrow\underline{E(n)\!}\,_{2r}$ and
an element $x\in H^*(\underline{E(n)\!}\,_{2r})$. It suffices to show that
$f^*(x)\in Ch_{E(n+1)}(X)$. Viewing $\underline{E(n)\!}\,_{2r}$ as the
homotopy
colimit as mentioned, compactness of $X$ means that the map $f$ factors
through
some intermediate space $\underline{BP\langle n\rangle\!}\,_{2r-k|v_n|}$ for
$k$ sufficiently large. Perhaps taking $k$ to be even larger, view this space
as a factor in
$\underline{BP\langle n+1\rangle\!}\,_{2r-k|v_n|}=
\underline{BP\langle n\rangle\!}\,_{2r-k|v_n|}\times
\underline{BP\langle n+1\rangle\!}\,_{2r-k|v_n|+|v_{n+1}|}$
by the Wilson splitting theorem. Then we have a factorisation of $f$ as
$$X\bra{f_1}\underline{BP\langle n+1\rangle\!}\,_{2r-k|v_n|}\bra{f_2}
\underline{E(n)\!}\,_{2r}$$
where $f_2$ is projection on the first factor followed by localisation with
respect to powers of $v_n$. Note that $f_2^*(x)$ will be an element of the
form $y\otimes 1$ in
$$H^*(\underline{BP\langle n+1\rangle\!}\,_{2r-k|v_n|})=
H^*(\underline{BP\langle n\rangle\!}\,_{2r-k|v_n|})\otimes
H^*(\underline{BP\langle n+1\rangle\!}\,_{2r-k|v_n|+|v_{n+1}|})\,.$$

Now consider the localisation of
$\underline{BP\langle n+1\rangle\!}\,_{2r-k|v_n|}$
with respect to powers of $v_{n+1}$. For $s$ sufficiently small, the map
$$\begin{array}{lcl}
\underline{BP\langle n+1\rangle\!}\,_{s} &=&
\underline{BP\langle n\rangle\!}\,_{s}
\times\underline{BP\langle n+1\rangle\!}\,_{s+|v_{n+1}|}\\&&\\
\qquad\bigg\downarrow v_{n+1}&&\\&&\\
\underline{BP\langle n+1\rangle\!}\,_{s-|v_{n+1}|} &=&
\underline{BP\langle n\rangle\!}\,_{s-|v_{n+1}|}
\times\underline{BP\langle n\rangle\!}\,_{s}
\times\underline{BP\langle n+1\rangle\!}\,_{s+|v_{n+1}|}
\end{array}$$
is the inclusion as the last two factors. In particular (again assuming $s$
sufficiently small), the $v_{n+1}$-localisation map
$\underline{BP\langle n+1\rangle\!}\,_{s}\lra\underline{E(n+1)\!}\,_{s}$
gives an epimorphism in cohomology. We thus have a diagram
$$\begin{array}{c}
X\bra{f_1}\underline{BP\langle n+1\rangle\!}\,_{2r-k|v_n|}\bra{f_2}
\underline{E(n)\!}\,_{2r}\ ,\\ \\
\!\!\!\!\!\!\!\!\!\!\!\!\bigg\downarrow g\\ \\
\underline{E(n+1)\!}\,_{2r-k|v_n|}
\end{array}$$
where $g$ is the $v_{n+1}$ localisation map, and an element
$z\in H^*(\underline{E(n+1)\!}\,_{2r-k|v_n|})$ with $g^*(z)=f_2^*(x)$.
Then $f^*(x)=(gf_1)^*(z)\in Ch_{E(n+1)}(X)$, showing $Ch_{E(n)}(X)$ to
be contained in $Ch_{E(n+1)}(X)$ as claimed.
\qed

\begin{rem}\label{six:rem}{\em
Recall the \lq\lq algebraic\rq\rq\ chromatic filtration of a finite spectrum
$X$ (for example, \cite[Def.\,7.5.3]{Rav}). This is given by a tower of
$p$-local spectra
$$L_0X\longleftarrow L_1X\longleftarrow
L_2X\longleftarrow\cdots\longleftarrow X$$
and is defined as the corresponding filtration of $\pi_*(X)$ by the kernels
of the homomorphisms $(\eta_X)_*\colon\pi_*(X)\rightarrow\pi_*(L_nX)$. Here
$L_n$ is localisation with respect to the spectra $E(n)$. The corresponding
tower in the category of {\em spaces\/} gives the filtration of $H^*(X)$ by
images of the $\eta_X^*$. Proposition~\ref{Bloc} shows the set of subalgebras
$\left\{Ch_{E(n)}(X)\right\}$ to be compatible with this filtration, but what
is not immediate for general spaces $X$ is that the subalgebras $Ch_{E(n)}(X)$
are nested. Theorem~\ref{six} shows that, under suitable finiteness
conditions,
there is a filtration
$$Ch_{E(1)}(X)\subset\cdots\subset Ch_{E(n)}(X)\subset
Ch_{E(n+1)}(X)\subset\cdots\subset H^*(X)$$
with a levelwise inclusion in the unstable algebraic chromatic filtration.
}\end{rem}

\begin{rem} {\em
In the light of the results of this section, the equivalences proved in
(\ref{fiveF}) and (\ref{fiveH}),
the inclusion of (\ref{lift}) and the non-equivalences in general of
$Ch_{E(n)}(-)$ and $Ch_{E(n+1)}(-)$
noted in \cite{CCC}, we are tempted to wonder whether, for suitable $E$ and
$X$, the subalgebra $Ch_E(X)$
in fact only depends on the Bousfield equivalence class of $E$. In the
first instance, \lq suitable\rq\
might include complex oriented $E$ and spaces $X$ which are either finite
CW complexes or
else are models of $BG$ for a finite group $G$}\end{rem}

\noindent{\small The Department of Mathematics and Computer Science,
University of
Leicester, University Road, Leicester, LE1 7RH, England.}

\noindent{\small Email: J.Hunton@mcs.le.ac.uk}

\vspace{0.2cm}\noindent{\small Department of Mathematics, University of
Wuppertal, Gau{\ss}stra{\ss}e~20, D-42097 Wuppertal, Germany.}

\noindent{\small Email: schuster@math.uni-wuppertal.de}

\begin{thebibliography}{99}

\itemsep=\smallskipamount
\newcommand{\article}[1]{#1}
\newcommand{\journal}[1]{\emph{#1}}


\bibitem{Adams}
J.~F. Adams.
\newblock \emph{Stable homotopy and generalised homology.}
\newblock Chicago Lectures in Math.\@, Chicago Univ.\@ Press, Chicago, 1974.

\bibitem{Atiyah} M.~F. Atiyah.
\newblock \article{Characters and the cohomology of finite groups.}
\newblock \journal{Inst.\ Hautes \'{E}tudes Sci.\ Publ.\ Math.}~\textbf{9}
(1961), 23--64.

\bibitem{AM}
M.~F. Atiyah and I.~G.~MacDonald.
\newblock \emph{Introduction to commutative algebra.}
\newblock Addison-Wesley, 1969.

\bibitem{bwa} A.~J. Baker and U. W\"urgler.
\newblock \article{Liftings of formal groups and the Artinian completion
of $v_n^{-1}BP$.}
\newblock \journal{Math.\ Proc.\ Camb.\ Phil.\ Soc.}~\textbf{106} (1989),
511--530.

\bibitem{bwb} A.~J. Baker and U. W\"urgler.
\newblock \article{Bockstein operations in Morava $K$-theories.}
\newblock \journal{Forum Math.}~\textbf{3} (1991), 543--560.

\bibitem{BH} M. Bendersky and J.~R. Hunton.
\newblock \article{On the coalgebraic ring and Bousfield-Kan spectral
sequence for a Landweber exact spectrum.}
\newblock In preparation.

\bibitem{bousfield} A. K. Bousfield.
\newblock \article{The localisation of spaces with respect to homology.}
\newblock \journal{Topology}~\textbf{14} (1975), 133--150.

\bibitem{carlson} Jon F. Carlson.
\newblock \article{Varieties and transfers.}
\newblock In: Proceedings of the Northwestern conference on cohomology of
groups (Evanston, Ill., 1985). \journal{J.\ Pure Appl.\ Algebra}~\textbf{44}
(1987), 99--105.

\bibitem{carlsson} G. Carlsson.
\newblock \article{Equivariant stable homotopy and Segal's Burnside ring
conjecture.}
\newblock \journal{Ann.\ of Math.}~\textbf{120} (1984), 189--224.

\bibitem{Evens} L. Evens.
\newblock \article{The cohomology ring of a finite group.}
\newblock \journal{Trans.\ Amer.\ Math.\ Soc.}~\textbf{101} (1961), 224--239.

\bibitem{VCh}
D.~J. Green and I.~J. Leary.
\newblock \article{The spectrum of the Chern subring.}
\newblock \journal{Comment.\@ Math.\@ Helv.\@}~\textbf{73} (1998), 406--426.

\bibitem{CCC} D.~J. Green, J. R. Hunton and B. Schuster.
\newblock \article{Chromatic characteristic classes in ordinary group
cohomology.}
\newblock Preprint at math.AT/0109019.

\bibitem{GLS} D.~J. Green, I. J. Leary and B. Schuster.
\newblock \article{The subring of group cohomology constructed by permutation
representations.}
\newblock To appear in \journal{Proc.\ Edinburgh Math.\ Soc.}

\bibitem{HopHov}
M.~J. Hopkins and M. Hovey.
\newblock \article{Spin cobordism determines real {$K$}-theory.}
\newblock \journal{Math.\ Z.}~\textbf{210} (1992), 181--196.

\bibitem{Hov} M. Hovey.
\newblock \article{Spin bordism and elliptic homology.}
\newblock \journal{Math.\ Z.}~\textbf{219} (1995), 163--170.

\bibitem{HH}
M.~J.~Hopkins and J.~R.~Hunton.
\newblock \article{On the structure of spaces representing a {L}andweber exact
cohomology theory.}
\newblock \journal{Topology}~\textbf{34} (1995), 29--36.

\bibitem{HKR}
M.~J. Hopkins, N.~J. Kuhn and D.~C. Ravenel.
\newblock \article{Generalized group characters and complex oriented cohomology
theories.}
\newblock \journal{J. Amer.\ Math.\ Soc.}~\textbf{13}
(2000), 553--594.

\bibitem{hun}
J.~R. Hunton.
\newblock \article{Complete cohomology theories and the homology of their
omega spectra.}
\newblock To appear in \journal{Topology}.

\bibitem{HR}
J.~R.~Hunton and N.~Ray.
\newblock \article{A rational approach to {H}opf rings.}
\newblock \journal{J.\ Pure Appl.\ Algebra}~\textbf{101}
(1995), 313--333.

\bibitem{HT2}
J.~R. Hunton and P.~R. Turner.
\newblock \article{The homology of spaces representing exact pairs of
homotopy functors.}
\newblock \journal{Topology}~\textbf{38} (1999), 621--634.

\bibitem{JW}
D. C. Johnson and W. S. Wilson.
\newblock \article{Projective dimension and Brown-Peterson homology.}
\newblock \journal{Topology}~\textbf{12} (1973), 327--353.

\bibitem{krst:hpt}
M. Kreck and S. Stolz.
\newblock \article{{$\mathbb{H}P^2$}-bundles and elliptic homology.}
\newblock \journal{Acta Math.}~\textbf{171} (1993), 231--261.

\bibitem{Landweber}
P.~S. Landweber.
\newblock \article{Cobordism and classifying spaces.}
\newblock In: A.~Liulevicius, Ed.
\newblock \emph{Algebraic Topology (Madison, 1970)\@.}
\newblock Proceedings of Symposia in Pure Mathematics\@ \textbf{22}, 125--129,
Amer.\@ Math.\@ Soc.\@, Providence, R.I., 1971.

\bibitem{nakaoka}
M. Nakaoka.
\newblock \article{Homology of the infinite symmetric group.}
\newblock \journal{Ann.\@ of Math.\@}~\textbf{73} (1961), 229--257.

\bibitem{Priddy}
S.~B.~Priddy.
\newblock \article{On $\Omega^\infty S^\infty$ and the infinite symmetric
group.}
\newblock In: A.~Liulevicius, Ed.
\newblock \emph{Algebraic Topology (Madison, 1970)\@.}
\newblock Proceedings of Symposia in Pure Mathematics~\textbf{22},
217--220, Amer.\@
Math.\@ Soc.\@, Providence, R.I., 1971.

\bibitem{Quillen}
D. Quillen.
\newblock \article{The spectrum of an equivariant cohomology ring I, II.}
\newblock \journal{Ann.\@ of Math.\@}~\textbf{94} (1971), 573--602.

\bibitem{Rav}
D.~C. Ravenel.
\newblock \emph{Nilpotence and periodicity in stable homotopy theory.}
\newblock Ann.\@ of Math.\@ Stud.~\textbf{128}, Princeton University
Press, 1992.

\bibitem{RW}
D.~C. Ravenel and W.~S. Wilson.
\newblock \article{The Hopf ring for complex cobordism.}
\newblock \journal{J.\@ Pure Appl.\@ Algebra}~\textbf{9} (1977), 241--280.

\bibitem{Thomas} C.~B. Thomas.
\newblock \emph{Characteristic classes and the cohomology of finite groups.}
\newblock Cambridge Studies in Advanced Mathematics, Cambridge University
Press, 1986.

\bibitem{Venkov}
B.~B.~Venkov.
\newblock \article{Cohomology algebras for some classifying spaces.}
\newblock \journal{Dokl.\ Akad.\ Nauk.\ SSSR}~\textbf{127} (1959),
943--944.

\bibitem{Venkov2}
B.~B.~Venkov.
\newblock \article{Characteristic classes for finite groups.}
\newblock \journal{Dokl.\ Akad.\ Nauk.\ SSSR}~\textbf{137} (1961), 1274--1277.


\bibitem{wilson2}
W.~S.~Wilson.
\newblock \article{The $\Omega$-spectrum for Brown-Peterson cohomology part
II.}
\newblock \journal{Amer.\@ J.\@ Math.}~\textbf{97} (1975),
101--123.

\bibitem{wkn}
W.~S. Wilson.
\newblock \article{The Hopf ring for Morava K-theory.}
\newblock \journal{Publ.\@ Res.\@ Inst.\@ Math.\@ Sci.,
Kyoto University}~\textbf{20} (1984), 1025--1036.

\bibitem{Yag:TAMS}
N. Yagita.
\newblock \article{Equivariant $BP$-cohomology for finite groups.}
\newblock \journal{Trans.\@ Amer.\@ Math.\@ Soc.}~\textbf{317} (1990),
485--499.

\bibitem{Yag:Land}
N. Yagita.
\newblock \article{The exact functor theorem for {$BP_*/I_n$}-theory.}
\newblock \journal{Proc.\ Japan Acad.}~\textbf{52} (1976), 1--3.


\end{thebibliography}
\end{document}